\documentclass[12pt]{article}
\usepackage{latexsym}
\usepackage{amsmath}
\usepackage{amssymb}
\usepackage{amsfonts}

\numberwithin{equation}{section}

\newtheorem{theorem}{Theorem}[section]
\newtheorem{lemma}[theorem]{Lemma}

\newtheorem{remark}[theorem]{Remark}
\newcommand{\eproof}{{\mbox{\ }~\hfill
\mbox{\large $\Box$} \par \vskip 10pt}}

\newcommand{\R}{{\mathbb R}}

\newcommand{\pf}{\noindent{\bf Proof}}
\renewcommand{\div}{{\rm div}}

\title{Optimal three-ball inequalities and quantitative uniqueness for the Lam\'e system with Lipschitz coefficients}
\author{Ching-Lung Lin\thanks{Department of Mathematics, National Cheng Kung University,
Tainan 701, Taiwan. Email:cllin2@mail.ncku.edu.tw}  \qquad Gen Nakamura\thanks{Department of
Mathematics, Hokkaido University, Sapporo 060-0810, Japan.
Partially supported by Grant-in-Aid for Scientific Research
(B)(2)(No. 14340038) of Japan Society for Promotion of Science.
(Email: gnaka@math.sci.hokudai.ac.jp)}\qquad Jenn-Nan
Wang\thanks{Department of Mathematics, Taida Institute of
Mathematical Sciences, NCTS (Taipei), National Taiwan University,
Taipei 106, Taiwan. Email:jnwang@math.ntu.edu.tw}}

\date{}

\begin{document}
\maketitle

\begin{abstract}
In this paper we study the local behavior of a solution to the
Lam\'e system with \emph{Lipschitz} coefficients in dimension $n\ge
2$. Our main result is the bound on the vanishing order of a
nontrivial solution, which immediately implies the strong unique
continuation property. This paper solves the open problem of the
strong uniqueness continuation property for the Lam\'e system with
Lipschitz coefficients in any dimension.
\end{abstract}

\section{Introduction}\label{sec1}
\setcounter{equation}{0}

Assume that $\Omega$ is a connected open set containing $0$ in
$\R^n$ for $n\geq 2$. Let $\lambda(x)$ and $\mu(x)$ be Lam\'e
coefficients in $C^{0,1}(\Omega)$ satisfying
\begin{equation}\label{1.1}
\begin{cases}
\mu(x)\geq\delta_0>0,\quad\quad\lambda(x)+2\mu(x)\geq\delta_0>0\quad\forall\
x\in\Omega,\\
\|\mu(x)\|_{C^{0,1}(\Omega)}+\|\lambda(x)\|_{C^{0,1}(\Omega)}\leq M_0.
\end{cases}
\end{equation}
The isotropic elasticity, which represents the displacement
equation of equilibrium, is given by
\begin{equation}\label{1.2}
\text{div}(\mu(\nabla u+(\nabla
u)^t))+\nabla(\lambda\text{div}u)=0\quad\text{in}\ \Omega,
\end{equation}
where $u=(u_1,u_2,\cdots,u_n)^t$ is the displacement vector and $(\nabla
u)_{jk}=\partial_ku_j$ for $j,k=1,2,\cdots,n$.

Results on the weak unique continuation for the Lam\'e system in
$\R^n,\ n\geq 2$, have been proved by Dehman and Robbiano for
$\lambda(x),\mu(x)\in C^{\infty}(\Omega)$ \cite{dero}, Ang, Ikehata,
Trong and Yamamoto for $\lambda\in C^2(\Omega), \mu(x)\in
C^3(\Omega)$ \cite{aity}, Weck for $\lambda(x),\mu(x)\in
C^2(\Omega)$ \cite{we1}, and Eller for $\lambda(x),\mu(x)\in
C^1(\Omega)$ \cite{el}.  As for the SUCP, it was proven by
Alessandrini and Morassi \cite{almo} in the case of
$\lambda(x),\mu(x)\in C^{1,1}(\Omega)$ and $n\ge 2$. Their proofs
were based on ideas developed by Garofalo and Lin \cite{gl1},
\cite{gl2}.  When Lam\'e coefficients are Lipschitz, i.e.,
$\lambda,\mu\in C^{0,1}(\Omega)$, the SUCP was established by the
first and third authors in \cite{lw05} for $n=2$. Later, the result
of \cite{lw05} was improved to $\mu\in C^{0,1}(\Omega)$ and
$\lambda$ being measurable by Escauriaza \cite{es}. In this work, we
completely resolve the SUCP problem for \eqref{1.2} when
$\lambda,\mu\in C^{0,1}(\Omega)$ and $n\ge 2$. It is important to
remark that in the three or higher dimensions, the Lipschitz
regularity assumption on the principal coefficients of the second
order elliptic equation is the minimal requirement for the unique
continuation property to hold \cite{pl}. Not only do we solve the
SUCP for the Lam\'e system with the minimal regularity assumption,
we also derive a quantitative form of the SUCP.

The ideas of our proof originate from our series papers on proving
quantitative uniqueness for elliptic equations or systems by the
method of Carleman estimates \cite{lin2}, \cite{lin3}, and
\cite{lin4}. In particular, the idea used in \cite{lin4} plays a key
role in our arguments here. Specifically, let us write \eqref{1.2}
into a non-divergence form:
\begin{equation}\label{uu}
\mu\Delta u +(\lambda+\mu)\nabla \div u+\nabla\lambda\ \div
u+(\nabla u+(\nabla u)^t)\nabla\mu=0.
\end{equation}
Letting $p=\div u$, taking divergence on \eqref{uu}, and using
\eqref{uu} for $\Delta u$, yields
\begin{eqnarray}\label{pp}
&&(\lambda+2\mu)\Delta
p+(\nabla\lambda-\mu^{-1}\lambda\nabla\mu)\cdot\nabla
p-(\mu^{-1}\nabla\mu\cdot\nabla\lambda)
\div u\notag\\
&&-\mu^{-1}\nabla\mu\cdot((\nabla u)+(\nabla
u)^t)\nabla\mu+\div\Big{(}\nabla\lambda\ \div u+(\nabla u+(\nabla
u)^t)\nabla\mu\Big{)}\notag\\
&&=0.
\end{eqnarray}
From \eqref{uu} and \eqref{pp}, we then obtain a system of equations
with the Laplacian as the principal part, namely,
\begin{equation}\label{up}
\begin{cases}
\Delta u+P_1(x,\partial)p+P_2(x,\partial)u=0,\\
\Delta p+Q_1(x,\partial)p+Q_2(x,\partial)u+\div G(x,u)=0,
\end{cases}
\end{equation}
where $P_j(x,\partial),Q_j(x,\partial)$, $j=1,2$, are first order
differential operators with at least essentially bounded
coefficients and
$$
G(x,u)=(\lambda+2\mu)^{-1}\Big{(}\nabla\lambda\ \div u+(\nabla
u+(\nabla u)^t)\nabla\mu\Big{)}.
$$
Note that system \eqref{up} is not decoupled.

To study the unique continuation problem for \eqref{1.2}, it
suffices to consider that for \eqref{up} with $p=\div u$. To tackle
this problem, we rely on suitable Carleman estimates. An important
step is to handle the second equation of \eqref{up}. The trick is to
use a Carleman estimate with the divergence operator on the right
hand side (see Lemma~\ref{lem2.4}). This idea was first introduced
in \cite{fale} and later used in \cite{Reg2} and \cite{lin4}. In
order to derive an upper bound on the vanishing order of a
nontrivial solution to \eqref{1.2}, it is also important to derive
\emph{optimal} three-ball inequalities.

We now state main results of the paper. Their proofs will be given
in the subsequent sections. Assume that there exists $0<R_0\le 1$
such that $B_{R_0}\subset\Omega$. Hereafter $B_r$ denotes an open
ball of radius $r>0$ centered at the origin.
\begin{theorem}\label{thm1.1}
There exists a positive number $\tilde{R}<1$, depending only on $n,M_0,\delta_0$,
such that if $\ 0<R_1<R_2<R_3\leq R_0$ and
$R_1/R_3<R_2/R_3<\tilde{R}$, then
\begin{equation}\label{1.3}
\int_{|x|<R_2}|u|^2dx\leq
{C}\left(\int_{|x|<R_1}|u|^2dx\right)^{\tau}\left(\int_{|x|<{R_3}}|u|^2dx\right)^{1-\tau}
\end{equation}
for $u\in H_{loc}^1({B}_{R_0})$ satisfying \eqref{1.2} in
${B}_{R_0}$, where the constant ${C}$ depends on $R_2/R_3$, $n$,
$M_0,\delta_0$, and $0<\tau<1$ depends on $R_1/R_3$, $R_2/R_3$,
$n,M_0,\delta_0$. Moreover, for fixed $R_2$ and $R_3$, the exponent
$\tau$ behaves like $1/(-\log R_1)$ when $R_1$ is sufficiently
small.
\end{theorem}
\begin{remark}\label{rem1.1}
We would like to emphasize that $C$ is independent of $R_1$ and
$\tau$ has the asymptotic $(-\log R_1)^{-1}$. These facts are
crucial in deriving an vanishing order of a nontrivial $u$ to
\eqref{1.1}. Due to the behavior of $\tau$, the three-ball
inequality is called optimal {\rm\cite{efv}}.
\end{remark}

\begin{theorem}\label{thm1.2}
Let $u\in H^1(\Omega)$ be a nontrivial solution of \eqref{1.2}, then
there exist positive constants $K$ and $m$, depending on
$n,M_0,\delta_0$ and $u$, such that
\begin{equation}\label{1.4}
\int_{|x|<R}|u|^2 dx\ge KR^m
\end{equation}
for all $R$ sufficiently small.
\end{theorem}
\begin{remark}\label{rem1.2}
Based on Theorem~\ref{thm1.1}, the constants $K$ and $m$ in
\eqref{1.4} are explicitly given by
$$
K=\int_{|x|<R_3}|u|^2dx
$$
and
$$
m=\tilde
C\log\Big{(}\frac{\int_{|x|<R_3}|u|^2dx}{\int_{|x|<R_2}|u|^2dx}\Big{)},
$$
where $\tilde C$ is a positive constant depending on $n,M_0,\delta_0$ and
$R_2/R_3$.
\end{remark}

\section{Carleman estimates}\label{sec2}
\setcounter{equation}{0}

In this section, we will derive two Carleman estimates. The first
one is taken from \cite{lin4}. Denote
$\varphi_{\beta}=\varphi_{\beta}(x) =\exp (-\beta\tilde{\psi}(x))$,
where $\beta>0$ and $\tilde{\psi}(x)=\log |x|+\log((\log |x|)^2)$.
Note that $\varphi_{\beta}$ is less singular than $|x|^{-\beta}$,
For simplicity, we denote $\psi(t)=t+\log t^2$, i.e.,
$\tilde{\psi}(x)=\psi(\log|x|)$. From now on, the notation
$X\lesssim Y$ or $X\gtrsim Y$ means that $X\le CY$ or $X\ge CY$ with
some constant $C$ depending only on $n$.

\begin{lemma}{\rm\cite[Lemma
2.1]{lin4}}\label{lem2.1}
There exist a sufficiently small $r_0>0$ depending on $n$ and a
sufficiently large $\beta_0>1$ depending on $n$ such that for all
$u\in U_{r_0}$ and $\beta\geq \beta_0$, we have that
\begin{equation}\label{2.5}
\beta  \int \varphi^2_\beta (\log|x|)^{-2}|x|^{-n}(|x|^{2}|\nabla
u|^2+|u|^2)dx\lesssim\int \varphi^2_\beta|x|^{-n}|x|^{4}|\Delta u|^2
dx,
\end{equation}
where $U_{r_0}=\{u\in C_0^{\infty}(\R^n\setminus\{0\}): \mbox{\rm
supp}(u)\subset B_{r_0}\}$.
\end{lemma}

To prove the second Carleman estimate, we need some preparations.
Firstly, we introduce polar coordinates in ${\mathbb R}^n \backslash
{\{0\}}$ by setting $x=r \omega$, with $r=|x|$,
$\omega=(\omega_1,\cdots,\omega_n)\in S^{n-1}$. Using new coordinate
$t=\log r$, we can see that
$$\frac{\partial}{\partial x_j}=e^{-t}(\omega_j \partial_t +\Omega_j),\quad 1\le j\le n,$$
where $\Omega_j$ is a vector field in $S^{n-1}$. We could check that
the vector fields $\Omega_j$ satisfy
$$\sum_{j=1}^n\omega_j\Omega_j=0\quad\text{and}\quad\sum_{j=1}^n\Omega_j\omega_j=n-1.$$
Since $r\rightarrow 0$ iff $t\rightarrow {-\infty}$, we are mainly
interested in values of $t$ near $-\infty$.

It is easy to see that
\begin{equation*}
\frac{\partial ^2}{\partial x_j \partial x_{\ell}}=e^{-2t}(\omega_j
\partial_t -\omega_j +\Omega_j)(\omega_{\ell} \partial_t +\Omega_{\ell}),\quad 1\le j,\ell\le n.
\end{equation*}
and, therefore, the Laplacian becomes
\begin{equation*}
e^{2t}\Delta =\partial^2_t +(n-2)\partial_t +\Delta_\omega,
\end{equation*}
where $\Delta_\omega=\Sigma^n_{j=1}\Omega^2_j$ denotes the
Laplace-Beltrami operator on $S^{n-1}$. We recall that the
eigenvalues of $-\Delta_\omega$ are $k(k+n-2), k\in \mathbb{N}$, and
the corresponding eigenspaces are $E_k$, where $E_k$ is the space of
spherical harmonics of degree $k$.  Let
$$
\Lambda=\sqrt{\frac{(n-2)^2}{4}-\Delta_{\omega}},
$$
then $\Lambda$ is an elliptic first-order positive
pseudodifferential operator in $L^2(S^{n-1})$. The eigenvalues of
$\Lambda$ are $k+\frac{n-2}{2}$ and the corresponding eigenspaces
are $E_k$ which represents the space of spherical harmonics of
degree $k$. Hence

\begin{equation}\label{2.6}
\Lambda=\Sigma_{k\geq 0}(k+\frac{n-2}{2})\pi_{k},
\end{equation}
where $\pi_k$ is the orthogonal projector on $E_k$. Denote
$$
L^{\pm}=\partial_t+\frac{n-2}{2}\pm\Lambda.
$$
Then it follows that
\begin{equation*}
e^{2t}\Delta =L^+L^-=L^-L^+.
\end{equation*}

We first recall a Carleman estimate proved in
\cite[Lemma~2.2]{lin4}.
\begin{lemma}\label{lem2.2}
There exists a sufficiently small number $t_0<0$ depending on $n$
such that for all $u\in V_{t_0}$, $\beta> 1$, we have that
\begin{equation}\label{2.7}
\sum_{j+|\alpha|\leq1}\beta^{1-2(j+|\alpha|)}\iint
t^{-2}\varphi^2_\beta|\partial_t^j\Omega^\alpha u|^2 dt d\omega
\lesssim \iint\varphi^2_\beta |L^- u|^2 dtd\omega,
\end{equation}
where $V_{t_0}=\{u(t,\omega)\in C_0^{\infty}((-\infty,t_0)\times
S^{n-1})\}$.
\end{lemma}

Next, we need an auxiliary Carleman estimate.
\begin{lemma}\label{lem2.3}
There exists a sufficiently small number $t_1<-2$ depending on $n$
such that for all $u\in V_{t_1}$, $g=(g_0,g_1,\cdots,g_n)\in
(V_{t_1})^{n+1}$ and $\beta> 1$, we have that
\begin{equation}\label{2.8}
\beta\iint \varphi^2_\beta |u|^2 dt d\omega \lesssim \iint
\varphi^2_\beta
(|L^+u+\partial_tg_0+\sum_{j=1}^n\Omega_jg_j|^2+\beta\|g\|^2)
dtd\omega.
\end{equation}
\end{lemma}
\pf. We shall prove this lemma following the lines of \cite[Lemma
2.2]{Reg2}. By defining $u=e^{\beta\psi(t)}v$ and
$g=e^{\beta\psi(t)}h$, \eqref{2.8} is equivalent to
\begin{equation}\label{2.9}
\beta\iint  |v|^2 dt d\omega \lesssim \iint(|L^+_\beta
v+(\partial_t+\beta\psi')h_0+\sum_{j=1}^n\Omega_jh_j|^2+\beta\|h\|^2)
dtd\omega,
\end{equation}
where $L^+_\beta v=e^{-\beta\psi(t)}L^+(e^{\beta\psi(t)} v)$ and
$h=(h_0,\cdots,h_n)$. By direct computations, we obtain that
\begin{eqnarray}\label{2.10}
&&\quad L^+_\beta v+(\partial_t+\beta\psi')h_0+\sum_{j=1}^n\Omega_jh_j\notag\\
&=&\partial_t v+\beta v+2\beta t^{-1}v+\frac{(n-2)}{2}v+\Lambda v+(\partial_t+\beta\psi')h_0+\sum_{j=1}^n\Omega_jh_j\notag\\
&=&T_\beta v+2\beta
t^{-1}v+(\partial_t+\beta\psi')h_0+\sum_{j=1}^n\Omega_jh_j
\end{eqnarray}
with $T_\beta v=\partial_t v+\beta v+\frac{(n-2)}{2}v+\Lambda v$.

For $v\in C^\infty_0(\mathbf{R}\times S^{n-1})$, we denote $\hat{v}$
its Fourier transformation with respect to $t$, then it follows from
\eqref{2.6} that
\begin{equation}\label{2.11}
T_\beta v(t,\omega)=(2\pi)^{-1}\sum_{k\geq 0}\int_{-\infty}^\infty
e^{i\sigma
t}(i\sigma+\beta+k+n-2)\pi_k\hat{v}(\sigma,\omega)d\sigma.
\end{equation}
It is easily seen that $T_\beta$ is invertible whose inverse is
given by
\begin{equation}\label{2.12}
T_\beta^{-1} v(t,\omega)=(2\pi)^{-1}\sum_{k\geq
0}\int_{-\infty}^\infty e^{i\sigma
t}(i\sigma+\beta+k+n-2)^{-1}\pi_k\hat{v}(\sigma,\omega)d\sigma.
\end{equation}

From \eqref{2.11}, \eqref{2.12} and Plancherel's theorem, we have
for $u\in C^\infty_0(\mathbf{R}\times S^{n-1})$ that
\begin{eqnarray}\label{2.13}
\begin{cases}
&\|T_\beta u\|_{L^2(\mathbf{R}\times S^{n-1})}\geq \beta\|u\|_{L^2(\mathbf{R}\times S^{n-1})},\\
&\beta\|T_\beta^{-1} u\|_{L^2(\mathbf{R}\times S^{n-1})}\leq \|u\|_{L^2(\mathbf{R}\times S^{n-1})},\\
&\|T_\beta^{-1}(\sum_{j+|\alpha|\leq1} \partial_t^j\Omega^\alpha
u)\|_{L^2(\mathbf{R}\times S^{n-1})}\lesssim
\|u\|_{L^2(\mathbf{R}\times S^{n-1})}.
\end{cases}
\end{eqnarray}
Combining \eqref{2.10} and \eqref{2.13}, we get that
\begin{eqnarray}\label{2.14}
&&\quad \|L^+_\beta v+(\partial_t+\beta\psi')h_0+\sum_{j=1}^n\Omega_jh_j\|_{L^2(\mathbf{R}\times S^{n-1})}\notag\\
&=&\|T_\beta v+2\beta t^{-1}v+(\partial_t+\beta\psi')h_0+\sum_{j=1}^n\Omega_jh_j\|_{L^2(\mathbf{R}\times S^{n-1})}\notag\\
&\ge&\beta\|v+T_\beta^{-1} \bigl(2\beta t^{-1}v+(\partial_t+\beta\psi')h_0+\sum_{j=1}^n\Omega_jh_j\bigr)\|_{L^2(\mathbf{R}\times S^{n-1})}\notag\\
&\geq& \beta\|v\|_{L^2(\mathbf{R}\times S^{n-1})}-\beta\|T_\beta^{-1} \bigl(2\beta t^{-1}v+(\partial_t+\beta\psi')h_0+\sum_{j=1}^n\Omega_jh_j\bigr)\|_{L^2(\mathbf{R}\times S^{n-1})}\notag\\
&\ge&\beta\|v\|_{L^2(\mathbf{R}\times S^{n-1})}-2\beta\|t^{-1}v\|_{L^2(\mathbf{R}\times S^{n-1})}\notag\\
&&-\beta\|T_\beta^{-1} \bigl((\partial_t+\beta\psi')h_0+\sum_{j=1}^n\Omega_jh_j\bigr)\|_{L^2(\mathbf{R}\times S^{n-1})}\notag\\
&\gtrsim& \beta\|v\|_{L^2(\mathbf{R}\times
S^{n-1})}-C\beta\|h\|_{L^2(\mathbf{R}\times S^{n-1})}.
\end{eqnarray}
In deriving \eqref{2.14}, we have used the fact that $v\in V_{t_1}$
with $t_1<-2$. Now dividing $\sqrt{\beta}$ on both sides of
\eqref{2.14} and squaring the new inequality, we have that
\begin{eqnarray}\label{2.15}
&&\quad\beta\iint  |v|^2 dt d\omega\notag\\
&\lesssim& \iint\beta^{-1}\bigl(\bigl |L^+_\beta v+(\partial_t+\beta\psi')h_0+\sum_{j=1}^n\Omega_jh_j\bigr |^2+\beta^2\|h\|^2\bigr ) dtd\omega\notag\\
&\lesssim& \iint\bigl(\bigl |L^+_\beta
v+(\partial_t+\beta\psi')h_0+\sum_{j=1}^n\Omega_jh_j\bigr
|^2+\beta\|h\|^2\bigr ) dtd\omega.
\end{eqnarray}
The proof is complete.\eproof

Now we are ready to prove our second Carleman estimate.
\begin{lemma}\label{lem2.4}
There exist a sufficiently small number $r_1>0$ depending on $n$ and
a sufficiently large number $\beta_1>3$ depending on $n$ such that
for all $w\in U_{r_1}$ and $f=(f_1,\cdots,f_n)\in (U_{r_1})^{n}$,
$\beta\geq \beta_1$, we have that
\begin{eqnarray}\label{2.16}
&&\int\varphi^2_\beta (\log|x|)^2(|x|^{4-n}|\nabla w|^2+|x|^{2-n}|w|^2)dx\notag\\
&\lesssim& \int \varphi^2_\beta
(\log|x|)^{4}|x|^{2-n}[(|x|^{2}\Delta w+|x|{\rm div}
f)^2+\beta\|f\|^2]dx,
\end{eqnarray}
where $U_{r_1}$ is defined as in Lemma~\ref{lem2.1}.
\end{lemma}

\pf. Replacing $\beta$ by $\beta+2$ in \eqref{2.16}, we see that it
suffices to prove
\begin{eqnarray}\label{2.17}
&&\int\varphi^2_\beta (\log|x|)^{-2}(|x|^{2}|\nabla w|^2+|w|^2)|x|^{-n}dx\notag\\
&\lesssim& \int \varphi^2_\beta [(|x|^{2}\Delta w+|x|{\rm div}
f)^2+\beta\| f\|^2 ] |x|^{-n}dx.
\end{eqnarray}
Working in polar coordinates and using the relation $e^{2t}\Delta
=L^+L^-$, \eqref{2.17} is equivalent to
\begin{eqnarray}\label{2.18}
&&\sum_{j+|\alpha|\leq1}\iint
t^{-2}\varphi^2_\beta|\partial_t^j\Omega^\alpha w|^2dtd\omega\notag\\
&\lesssim& \iint \varphi^2_\beta (|L^+L^-w+\partial_t(\sum_{j=1}^n
\omega_jf_j)+\sum_{j=1}^n\Omega_jf_j|^2+\beta\|f\|^2)dtd\omega.
\end{eqnarray}
Applying Lemma \ref{lem2.3} to $u=L^-w$ and $g=(\sum_{j=1}^n
\omega_jf_j,f_1,\cdots,f_n)$ yields
\begin{eqnarray}\label{2.19}
&&\beta\iint
\varphi^2_\beta |L^-w|^2 dt d\omega\notag\\
&\lesssim& \iint\varphi^2_\beta
(|L^+L^-w+\partial_t(\sum_{j=1}^n\omega_jf_j)+\sum_{j=1}^n\Omega_jf_j|^2+\beta\|f\|^2)
dtd\omega.
\end{eqnarray}
Now \eqref{2.18} is an easy consequence of \eqref{2.7} and
\eqref{2.19}.\eproof

\section{Proof of Theorem \ref{thm1.1} and Theorem \ref{thm1.2}}\label{sec4}
\setcounter{equation}{0}

This section is devoted to the proofs of Theorem~\ref{thm1.1} and
\ref{thm1.2}. To begin, we state an interior estimate for the Lam\'e
system \eqref{1.2}. For fixed $a_3<a_1<a_2<a_4$, there exists a
constant $\tilde C$ such that
\begin{equation}\label{2.3}
\int_{a_1r<|x|<a_2r}||x|^{|\alpha|}D^{\alpha} u|^2dx\le \tilde
C\int_{a_3r<|x|<a_4r}|u|^2dx,\quad|\alpha|\le 2
\end{equation}
for all sufficiently small $r$. Estimate \eqref{2.3} can be proved
by repeating the arguments of Corollary~17.1.4 in \cite{Hor3}.
Plugging $p=\div u$ in \eqref{2.3} yields
\begin{equation}\label{2.4}
\int_{a_1r<|x|<a_2r}||x|^{|\alpha|+1}D^{\alpha} p|^2dx\le \tilde
C\int_{a_3r<|x|<a_4r}|u|^2dx,\quad|\alpha|\le 1.
\end{equation}

To proceed the proof, let us first consider the case where
$0<R_1<R_2<R<R_0$. The constant $R$ will be determined later. Since
$u\in H_{loc}^1(B_{R_0})$, the elliptic regularity theorem for
\eqref{1.2} implies $u\in H^2_{loc}(B_{R_0})$. Therefore, to use
estimate \eqref{2.5}, we simply cut-off $u$. So let $\chi(x)\in
C^{\infty}_0 ({\mathbb R}^n)$ satisfy $0\le\chi(x)\leq 1$ and
\begin{equation*}
\chi (x)=
\begin{cases}
\begin{array}{l}
0,\quad |x|\leq R_1/e,\\
1,\quad R_1/2<|x|<eR_2,\\
0,\quad |x|\geq 3R_2,
\end{array}
\end{cases}
\end{equation*}
where $e=\exp(1)$. We first choose a small $R$ such that
$R\le\min\{r_0,r_1\}/3=\tilde R_0$, where $r_0$ and $r_1$ are
constants appeared in \eqref{2.5} and \eqref{2.16}. Hence $\tilde
R_0$ depends on $n$. It is easy to see that for any multiindex
$\alpha$
\begin{equation}\label{3.1}
\begin{cases}
|D^{\alpha}\chi|=O(R_1^{-|\alpha|})\ \text{for all}\ R_1/e\le |x|\le R_1/2\\
|D^{\alpha}\chi|=O(R_2^{-|\alpha|})\ \text{for all}\ eR_2\le |x|\le
3R_2.
\end{cases}
\end{equation}
Applying \eqref{2.5} to $\chi u$ gives
\begin{equation}\label{3.2}
{C}_1\beta  \int (\log|x|)^{-2}\varphi^2_\beta
|x|^{-n}(|x|^{2}|\nabla (\chi u)|^2+|\chi u|^2)dx \leq \int
\varphi^2_\beta |x|^{-n}|x|^{4}|\Delta (\chi u)|^2 dx.
\end{equation}
From now on, $C_1,C_2,\cdots$ denote general constants whose
dependence will be specified whenever necessary. Next we want to
apply \eqref{2.16} to $w=\chi p$ and $f=|x|\chi G$. Since $u\in
H^2_{loc}$ and $p=\div u\in H^1_{loc}$, in view of the second
equation of \eqref{up}, by the standard limiting argument,
\eqref{2.16} holds true for $(w,f)$ above. Thus, we get that
\begin{eqnarray}\label{3.3}
&&{C}_2\int\varphi^2_\beta (\log|x|)^2(|x|^{4-n}|\nabla (\chi p)|^2+|x|^{2-n}|\chi p|^2)dx\notag\\
&\leq& \int \varphi^2_\beta(\log|x|)^{4}|x|^{2-n}[|x|^{2}\Delta(\chi p)+|x|{\rm div} (|x|\chi G)]^2dx\notag\\
&&+\beta\int \varphi^2_\beta (\log|x|)^{4}|x|^{2-n}\| |x|\chi G\|^2
dx.
\end{eqnarray}
Combining \eqref{3.2} and \eqref{3.3}, we obtain that
\begin{eqnarray}\label{3.4}
&&\beta \int_{R_1/2<|x|<eR_2}(\log|x|)^{-2}\varphi^2_\beta |x|^{-n}(|x|^{2}|\nabla u|^2+| u|^2)dx\notag\\
&&+\int_{R_1/2<|x|<eR_2}(\log|x|)^2\varphi^2_\beta |x|^{-n}(|x|^{4}|\nabla p|^2+|x|^{2}|p|^2)dx\notag\\
&\leq&\beta \int\varphi^2_\beta (\log|x|)^{-2}|x|^{-n}(|x|^{2}\nabla (\chi u)|^2+|\chi u|^2)dx\notag\\
&&+ \int(\log|x|)^2\varphi^2_\beta |x|^{-n}(|x|^{4}|\nabla (\chi p)|^2+|x|^{2}|\chi p|^2)dx\notag\\
&\leq &C_{3}\int\varphi^2_\beta |x|^{-n}|x|^{4}|\Delta (\chi u)|^2 dx\notag\\
&&+ C_{3}\int(\log|x|)^{4}\varphi^2_\beta |x|^{-n}[|x|^{3}\Delta (\chi p)+|x|^2{\rm div} (|x|\chi G)]^2dx\notag\\
&&+\beta C_{3}\int(\log|x|)^{4}\varphi^2_\beta |x|^{-n}\| |x|^2\chi
G\|^2 dx.
\end{eqnarray}
By \eqref{1.1}, \eqref{up}, and estimates \eqref{3.1},
we deduce from \eqref{3.4} that
\begin{eqnarray}\label{3.5}
&&\beta \int_{R_1/2<|x|<eR_2}
(\log|x|)^{-2}\varphi^2_\beta |x|^{-n}(|x|^{2}|\nabla u|^2+|u|^2)dx\notag\\
&&+\int_{R_1/2<|x|<eR_2}(\log|x|)^{2}\varphi^2_\beta |x|^{-n}(|x|^{4}|\nabla p|^2+|x|^{2}|p|^2)dx\notag\\
&\leq& C_4\int_{R_1/2<|x|<eR_2}\varphi^2_\beta |x|^{-n}|x|^{4}(|\nabla u|^2+|\nabla p|^2) dx\notag\\
&&+C_{4}\int_{R_1/2<|x|<eR_2}(\log|x|)^{4}|x|^2\varphi^2_\beta |x|^{-n}(|x|^{4}|\nabla u|^2+|x|^{4}|\nabla p|^2)dx\notag\\
&&+C_{4}\beta \int_{R_1/2<|x|<eR_2}(\log|x|)^{4}|x|^2\varphi^2_\beta |x|^{-n}|x|^{2}|\nabla u|^2dx\notag\\
&&+C_4\int_{\{R_1/e\le |x|\le R_1/2\}\cup\{eR_2\le |x|\le 3R_2\}}(\log|x|)^{4}\varphi^2_\beta |x|^{-n}|\tilde{U}(x)|^2 dx\notag\\
&&+C_{4}\beta \int_{\{R_1/e\le |x|\le R_1/2\}\cup\{eR_2\le
|x|\le3R_2\}}(\log|x|)^{4}\varphi^2_\beta |x|^{-n}|\tilde{U}(x)|^2 dx,
\end{eqnarray}
where $|\tilde{U}(x)|^2=|x|^{4}|\nabla
p|^2+|x|^{2}|p|^2+|x|^{2}|\nabla u|^2+|u|^2$ and the positive
constant $C_4$ only depends on $n,M_0,\delta_0$.

Now letting $R$ small enough, say $R<\tilde R_1$, such that $2C_4(\log(eR))^6(eR)^{2}\leq
1$ and $(\log(eR))^2\geq 2C_4$, then the first three terms on the right hand side of \eqref{3.5}
can be absorbed by the left hand side of \eqref{3.5}. Also, it is
easy to check that there exists $\tilde R_2>0$, depending on $n$,
such that for all $\beta>0$, both
$(\log|x|)^{-2}|x|^{-n}\varphi_{\beta}^2(|x|)$ and
$(\log|x|)^{4}|x|^{-n}\varphi_{\beta}^2(|x|)$ are decreasing
functions in $0<|x|<\tilde R_2$. So we choose a small $R<\tilde
R_3$, where $\tilde R_3=\min\{\tilde R_2/3,\tilde R_1,\tilde R_0\}$. It is
clear that $\tilde R_3$ depends on $n,M_0,\delta_0$. With the
choices described above, we obtain from \eqref{3.5} that

\begin{eqnarray}\label{3.6}
&&R_2^{-n}(\log R_2)^{-2}\varphi^2_\beta(R_2)\int_{R_1/2<|x|<R_2}|u|^2dx\notag\\
&\leq &\int_{R_1/2<|x|<eR_2}(\log|x|)^{-2}\varphi^2_\beta |x|^{-n}|u|^2dx\notag\\
&\leq &C_5 \int_{\{R_1/e\le |x|\le R_1/2\}\cup\{eR_2\le |x|\le 3R_2\}}(\log|x|)^{4}\varphi^2_\beta |x|^{-n}|\tilde{U}|^2 dx\notag\\
&\leq &C_5 (\log(R_1/e))^{4}(R_1/e)^{-n}\varphi^2_\beta(R_1/e)\int_{\{R_1/e\le |x|\le R_1/2\}}|\tilde{U}|^2 dx\notag\\
&&+C_5(\log(eR_2))^{4}(eR_2)^{-n}\varphi^2_\beta(eR_2)\int_{\{eR_2\le|x|\le
3R_2\}}|\tilde{U}|^2 dx.
\end{eqnarray}

It follows from \eqref{2.3} and \eqref{2.4} that
\begin{eqnarray}\label{3.7}
&&R_2^{-2\beta-n}(\log R_2)^{-4\beta-2}\int_{R_1/2<|x|<R_2}|u|^2dx\notag\\
&\leq & C_{6}(\log(R_1/e))^{4}(R_1/e)^{-n}\varphi^2_\beta(R_1/e)\int_{\{R_1/4\le |x|\le R_1\}}|u|^2 dx\notag\\
&&+C_{6} (\log(eR_2))^{4}(eR_2)^{-n}\varphi^2_\beta(eR_2)\int_{\{2R_2\le |x|\le 4R_2\}}|u|^2 dx\notag\\
&= & C_{6}(\log(R_1/e))^{-4\beta+4}(R_1/e)^{-2\beta-n}\int_{\{R_1/4\le |x|\le R_1\}}|u|^2 dx\notag\\
&&+C_{6}(\log(eR_2))^{-4\beta+4}(eR_2)^{-2\beta-n}\int_{\{2R_2\le
|x|\le 4R_2\}}|u|^2 dx.
\end{eqnarray}
Replacing $2\beta+n$ by $\beta$, \eqref{3.7} becomes
\begin{eqnarray}\label{3.8}
&&R_2^{-\beta}(\log R_2)^{-2\beta+2n-2}\int_{R_1/2<|x|<R_2}|u|^2dx\notag\\
&\leq & C_{7}(\log(R_1/e))^{-2\beta+2n+4}(R_1/e)^{-\beta}\int_{\{R_1/4\le |x|\le R_1\}}|u|^2 dx\notag\\
&&+C_{7}(\log(eR_2))^{-2\beta+2n+4}(eR_2)^{-\beta}\int_{\{2R_2\le
|x|\le4R_2\}}|u|^2 dx.
\end{eqnarray}
Dividing $R_2^{-\beta}(\log R_2)^{-2\beta+2n-2}$ on the both sides
of \eqref{3.8} and if $\beta\geq n+2$, we have that
\begin{eqnarray}\label{3.9}
&&\int_{R_1/2<|x|<R_2}|u|^2dx\notag\\
&\leq & C_{8}(\log R_2)^6(eR_2/R_1)^{\beta}\int_{\{R_1/4\le |x|\le R_1\}}|u|^2 dx\notag\\
&&+C_{8}(\log R_2)^{6}(1/e)^\beta [(\log R_2/\log(eR_2))^{2}]^{\beta-n-2}\int_{\{2R_2\le |x|\le 4R_2\}}|u|^2 dx\notag\\
&\leq & C_{8}(\log R_2)^{6}(eR_2/R_1)^{\beta}\int_{\{R_1/4\le |x|\le R_1\}}|u|^2 dx\notag\\
&&+C_{8}(\log R_2)^{6}(4/5)^\beta \int_{\{2R_2\le |x|\le4R_2\}}|u|^2
dx.
\end{eqnarray}
In deriving the second inequality above, we use the fact that
$$
\frac{\log R_2}{\log(eR_2)}\to 1\quad\text{as}\quad R_2\to 0,
$$
and thus
$$
\frac{1}{e}\cdot\frac{\log R_2}{\log(eR_2)}<\frac{4}{5}
$$
for all $R_2<\tilde R_4$, where $\tilde R_4$ is sufficiently small.
We now take $\tilde R=\min\{\tilde R_3,\tilde R_4\}$, which depends
on $n,M_0,\delta_0$.

Adding $\int_{|x|<{R_1/2}} |u|^2 dx$ to both sides of \eqref{3.9}
leads to
\begin{eqnarray}\label{3.10}
\int_{|x|<R_2}|u|^2dx&\leq&  C_{9}(\log R_2)^{6}(eR_2/R_1)^{\beta}\int_{|x|\le R_1}|u|^2dx\notag\\
&&+C_{9}(\log R_2)^{6}(4/5)^\beta \int_{|x|\le 1}|u|^2 dx.
\end{eqnarray}
It should be noted that \eqref{3.10} holds for all
$\beta\ge\tilde\beta$ with $\tilde\beta$ depending only on
$n,M_0,\delta_0$. For simplicity, by denoting
\begin{equation*}
E(R_1,R_2)=\log(eR_2/R_1),\quad B=\log (5/4),
\end{equation*}
\eqref{3.10} becomes
\begin{eqnarray}\label{3.11}
&&\int_{|x|<R_2}|u|^2dx\notag\\
&\leq & C_{9}(\log R_2)^{6}\Big{\{}\exp(E\beta)\int_{|x|<{R_1}}|u|^2 dx+\exp(-B\beta) \int_{|x|<1} |u|^2 dx\Big{\}}.\notag\\
\end{eqnarray}

To further simplify the terms on the right hand side of
\eqref{3.11}, we consider two cases. If $\int_{|x|<{R_1}} |u|^2
dx\ne 0$ and
$$\exp{(E\tilde\beta)}\int_{|x|<{R_1}} |u|^2 dx<\exp{(-B\tilde\beta)}\int_{|x|<{1}} |u|^2 dx,$$
then we can pick a $\beta>\tilde\beta$ such that
$$
\exp{(E\beta)}\int_{|x|<{R_1}} |u|^2
dx=\exp{(-B\beta)}\int_{|x|<{1}} |u|^2 dx.
$$
Using such $\beta$, we obtain from \eqref{3.11} that
\begin{eqnarray}\label{3.12}
&&\int_{|x|<R_2}|u|^2dx\notag\\
&\leq& 2C_{9}(\log R_2)^{6}\exp{(E\beta)} \int_{|x|<{R_1}} |u|^2dx\notag\\
&=& 2C_{9}(\log
R_2)^{6}\left(\int_{|x|<{R_1}}|u|^2dx\right)^{\frac{B}{E+B}}\left(\int_{|x|<{1}}|u|^2dx\right)^{\frac{E}{E+B}}.
\end{eqnarray}
If $\int_{|x|<{R_1}} |u|^2 dx= 0$, then letting $\beta\to\infty$ in
\eqref{3.11} we have $\int_{|x|<R_2}|u|^2dx=0$ as well. The
three-ball inequality obviously holds.

On the other hand, if
$$ \exp{(-B\tilde\beta)}\int_{|x|<{1}} |u|^2dx\leq\exp{(E\tilde\beta)}\int_{|x|<{R_1}} |u|^2 dx,$$
then we have
\begin{eqnarray}\label{3.13}
&&\int_{|x|<{R_2}}|u|^2 dx\notag\\
&\leq & \left(\int_{|x|<1}|u|^2dx\right)^{\frac{B}{E+B}}\left(\int_{|x|<1}|u|^2dx\right)^{\frac{E}{E+B}}\notag\\
&\leq &
\exp{(B\tilde\beta)}\left(\int_{|x|<{R_1}}|u|^2dx\right)^{\frac{B}{E+B}}\left(\int_{|x|<1}|u|^2dx\right)^{\frac{E}{E+B}}.
\end{eqnarray}
Putting together \eqref{3.12}, \eqref{3.13}, and setting
$C_{10}=\max\{2C_{9}(\log R_2)^{6},\exp{(\tilde\beta\log(5/4))}\}$,
we arrive at
\begin{equation}\label{3.14}
\int_{|x|<{R_2}}|u|^2 dx \le
C_{10}\left(\int_{|x|<{R_1}}|u|^2dx\right)^{\frac{B}{E+B}}\left(\int_{|x|<1}|u|^2dx\right)^{\frac{E}{E+B}}.
\end{equation}
It is readily seen that $\frac{B}{E+B}\approx (\log (1/R_1))^{-1}$
when $R_1$ tends to $0$.

Now for the general case, we consider $0<R_1<R_2<R_3<1$ with
$R_1/R_3<R_2/R_3\le \tilde R$, where $\tilde R$ is given as above.
By scaling, i.e. defining $\widehat{u}(y):=u(R_3y)$,
$\widehat{\lambda}(y):=\lambda(R_3y)$ and
$\widehat{\mu}(y)=\mu(R_3y)$,  \eqref{3.14}  becomes
\begin{equation}\label{3.15}
\int_{|y|<{R_2/R_3}}|\widehat{u}(y)|^2 dy \leq
C_{11}(\int_{|y|<{R_1/R_3}}|\widehat{u}(y)|^2dy)^{\tau}(\int_{|y|<1}|\widehat{u}(y)|^2dy)^{1-\tau},
\end{equation}
where $$\tau=B/[E(R_1/R_3,R_2/R_3)+B],$$ $$C_{11}=\max\{2C_{9}(\log
R_2/R_3)^{6},\exp{(\tilde\beta\log(5/4))}\}.$$ Note that $C_{11}$ is
independent of $R_1$. Restoring the variable $x=R_3y$ in
\eqref{3.15} gives
\begin{equation*}
\int_{|x|<{R_2}}|u|^2 dx \leq
C_{11}(\int_{|x|<{R_1}}|u|^2dx)^{\tau}(\int_{|x|<{R_3}}|u|^2dx)^{1-\tau}.
\end{equation*}
The proof of Theorem \ref{thm1.1} is complete.

We now turn to the proof of Theorem \ref{thm1.2}. We fix $R_2$,
$R_3$ in Theorem~\ref{thm1.1} and  define
$$
\widetilde{u}(x):=u(x)/\sqrt{\int_{|x|<{R_2}}|u|^2 dx}.
$$
Note that $\int_{|x|<{R_2}}|\widetilde{u}|^2 dx=1$. From the
three-ball inequality \eqref{1.2}, we have that
\begin{equation}\label{3.17}
1 \leq
C(\int_{|x|<{R_1}}|\widetilde{u}|^2dx)^{\tau}(\int_{|x|<{R_3}}|\widetilde{u}|^2dx)^{1-\tau}.
\end{equation}
Raising both sides by $1/\tau$ yields that
\begin{equation}\label{3.18}
\int_{|x|<{R_3}}|\widetilde{u}|^2dx \leq
(\int_{|x|<{R_1}}|\widetilde{u}|^2dx)(C\int_{|x|<{R_3}}|\widetilde{u}|^2dx)^{1/\tau}.
\end{equation}
In view of the formula for $\tau$, we can deduce from \eqref{3.18}
that
\begin{equation}\label{3.19}
\int_{|x|<{R_3}}|\widetilde{u}|^2dx \leq
(\int_{|x|<{R_1}}|\widetilde{u}|^2dx)(1/R_1)^{\tilde
C\log(\int_{|x|<{R_3}}|\widetilde{u}|^2dx)},
\end{equation}
where $\tilde C$ is a positive constant depending on $n,M_0,\delta_0$ and
$R_2/R_3$. Consequently, \eqref{3.19} is equivalent to
$$
(\int_{|x|<R_3}|u|^2dx) R_1^m\le \int_{|x|<R_1}|u|^2dx
$$
for all $R_1$ sufficiently small, where
$$
m=\tilde
C\log\Big{(}\frac{\int_{|x|<R_3}|U|^2dx}{\int_{|x|<R_2}|U|^2dx}\Big{)}.
$$
We now end the proof of Theorem \ref{thm1.2}.

\section*{Acknowledgements}
The first and third authors were supported in part by the National
Science Council of Taiwan.


\begin{thebibliography}{50}
\bibitem{almo}
 G. Alessandrini and A. Morassi, \emph{Strong unique continuation for the
 Lam\'e
 system of elasticity}, Comm. in PDE., {\bf 26}, 1787-1810, 2001.
\bibitem{aity}
   D.D. Ang, M. Ikehata, D.D. Trong and M. Yamamoto, \emph{Unique
   continuation for a stationary isotropic Lam\'e system with
   varaiable coefficients}, Comm. in PDE, {\bf 23}, 371-385, 1998.
\bibitem{dero}
 B. Dehman and L. Robbiano, \emph{La propri\'{e}t\'{e} du prolongement unique pour
 un syst\`{e}me elliptique: le syst\`{e}me Lam\'{e}}, J. Math. Pures Appl.,
 {\bf 72}, 475-492, 1993.
\bibitem{gl1}
N. Garofalo and F.H. Lin, \emph{Monotonicity properties of
variational integrals, $A\sb p$ weights and unique continuation},
Indiana Univ. Math. J., {\bf 35}, 245-268, 1986.
\bibitem{gl2}
N. Garofalo and F.H. Lin  \emph{Unique continuation for elliptic
operators: a geometric-variational approach}, Comm. Pure Appl.
Math., {\bf 40}, 347-366, 1987.



\bibitem{el}
M. Eller, \emph{Carleman estimates for some elliptic systems},
Journal of Physics: Conference Series \textbf{124} (2008) 012023.

\bibitem{es}
L. Escauriaza, \emph{ Unique continuation for the system of elasticity in the plane}, Proc. Amer. Math. Soc., \textbf{134} (2006), 2015-2018.

\bibitem{efv}
L. Escauriaza, F.J. Fern\'andez, and S. Vessella, \emph{Doubling
properties of caloric functions}, Appl. Anal., \textbf{85} (2006),
205-223.
\bibitem{fale}
C. Fabre and G. Lebeau, \emph{Prolongement unique des solutions de
l'\'{e}quation de Stokes}, Comm. in PDE, \textbf{21} (1996),
573-596.

\bibitem{Hor3}
L. H\"{o}rmander, {"The analysis of linear partial differential
operators"}, Vol. 3, Springer-Verlag, Berlin/New York, 1985.
\bibitem{lw05}
C.L. Lin and J.N. Wang, \emph{Strong unique continuation for the
Lam\'e system with Lipschitz coefficients}, Math. Ann., \textbf{331}
(2005), 611-629.
\bibitem{lin2}
C.L. Lin, G. Nakamura and J.N. Wang \emph{Quantitative uniqueness
for second order elliptic operators with strongly singular
coefficients}, Preprint(2008), http://arxiv.org/abs/0802.1983.
\bibitem{lin3}
C.L. Lin, S. Nagayasu and J.N. Wang  \emph{Quantitative uniqueness
for the power of Laplacian with singular coefficients}, Preprint(2008),
http://arxiv.org/abs/0803.1012.
\bibitem{lin4}
C.L. Lin and J.N. Wang  \emph{Optimal three-ball inequalities and quantitative uniqueness for the Stokes system},
Preprint(2008), http://arxiv.org/abs/0812.3730.
\bibitem{pl}
A. Plis, \emph{On non-uniqueness in Cauchy problem for an elliptic
second order differential equation}, Bull. Acad. Polon. Sci. Ser.
Sci. Math. Astronom. Phys., \textbf{11} (1963), 95-100.
\bibitem{Reg2}
R. Regbaoui, {\it Strong unique continuation for stokes equations},
Comm. in PDE {\bf 24} (1999), 1891-1902.
\bibitem{we1}
   N. Weck, \emph{Unique continuation for systems with Lam\'{e} principal part},
   Math. Methods Appl. Sci., {\bf 24}, 595--605, 2001.


\end{thebibliography}
\end{document}